\documentclass[11pt]{article}
\usepackage{amssymb}
\usepackage{amsthm}
\newtheorem{Theorem}{Theorem}

\newtheorem{Question}{Question}

\title{Towards bounded negativity of self-intersection\\
on general blown-up projective planes}
\author{Claudio Fontanari}
\date{}

\begin{document}
\maketitle

\begin{abstract}
We address the problem of bounding from below the self-intersection 
of integral curves on the projective plane blown-up at general points. 
In particular, by applying classical deformation theory we obtain 
the expected bound in the case of either high ramification or low multiplicity.
\end{abstract}

\section{Introduction}

Let $S$ be the blow-up of the complex projective plane $\mathbb{P}^2$ 
at $n \ge 1$ general points $p_1, \ldots, p_n$. Denote by $H$ the hyperplane 
class and by $E_i$ the exceptional divisor for $i=1, \ldots, n$.

The following natural problem seems to be still widely open:

\begin{Question} \label{conjecture}
Is there a constant $c_n$ depending only on $n$ such that 
the self-intersection number $C^2$ satisfies $C^2 \ge c_n$ 
for every integral curve $C \subset S$?
\end{Question}

Renewed interest in this question has been recently witnessed 
by both Joe Harris and Brian Harbourne (see \cite{H}, 
Question on p. 24, and \cite{H2}, Conjectures I.2.1 and I.2.7). 

According to \cite{H1}, the celebrated Segre-Harbourne-Gimigliano-Hirschowitz 
(SHGH) Conjecture (see for instance \cite{H1}, Conjecture 3.1) turns out to be 
equivalent to the sharp inequality $C^2 \ge g -1$, where $g$ is the arithmetic 
genus of $C$, hence the expected lower bound is precisely $C^2 \ge -1$.

It is indeed well-known that $C^2 \ge -1$ if $C$ is rational (see for instance 
\cite{dF1}, Proposition 2.4). The main result of \cite{dF1} 
shows in particular that $C^2 \ge -1$ if 
$C \in \vert dH - \sum_{i=1}^n m_i E_i \vert$ with $m_k = 2$ for some $k$. 
The Mori-theoretic point of view introduced in \cite{dF1} has been further 
developed in \cite{M} and \cite{dF2}.

Here we generalize \cite{dF1}, Theorem 2.5, in two different directions:

\begin{Theorem}\label{one}
Let $\Gamma$ be an integral curve in $\mathbb{P}^2$ and let 
$C$ be its strict transform. 
If $\Gamma$ has at most two tangent directions at $p_k$ for some $k$, 
then $C^2 \ge -1$.
\end{Theorem}

\begin{Theorem}\label{two}
Let $\Gamma$ be an integral curve in $\mathbb{P}^2$ and let 
$C \in \vert dH - \sum_{i=1}^n m_i E_i \vert$ be its strict transform. 
If $m_k \le 3$ for some $k$, then $C^2 \ge -1$.
\end{Theorem}

Our main tool is classical deformation theory. In particular, we follow 
the well-established approach of Xu (see \cite{Xu1}, \cite{Xu2}). 
We also exploit some recent refinements which appeared in \cite{B}, Lemma 3, 
and \cite{K}, Theorem A. For further applications of deformation theory 
to linear systems of divisors we refer to \cite{CM} and \cite{CC}, \S 2.    

We work over the complex field $\mathbb{C}$.

We are grateful to Edoardo Ballico, Ciro Ciliberto and Edoardo Sernesi 
for stimulating discussions on these topics and to the anonymous referee 
for her/his valuable remarks.

This research has been partially supported by GNSAGA of INdAM and MIUR Cofin 2008 - 
Geo\-metria delle variet\`{a} algebriche e dei loro spazi di moduli (Italy).

\section{The proofs}

\emph{Proof of Theorem \ref{one}.} Let $\Gamma \in 
\vert dH - (\sum_{i \ne k} m_i p_i - (m_k -1) p_k) \vert$
for appropriate choices of $d$ and $m_j$. The proof of \cite{Xu1}, Lemma 1 
(see also \cite{EL}, Lemma 1.1) shows that the linear system 
$\vert dH - (\sum_{i \ne k} m_i p_i - (m_k -1) p_k) \vert$ contains 
a curve $\Gamma' \ne \Gamma$. More explicitly, if $\Gamma = \{ F(X,Y,Z)=0 \}$,
$p_k(t):=[a(t),b(t),1]$ with $p_k(0)=p_k$ and $p_i(t):=p_i$ 
for every $i \ne k$, then there is a deformation $\Gamma_t = \{ F_t(X,Y,Z) = 0 \}$ 
of $\Gamma$ such that $F_0(X,Y,Z)=F(X,Y,Z)$ and $\Gamma_t$ passes through $p_i(t)$ 
with multiplicity $m_i$ for every $i=1, \ldots, n$ and every $t$
in a neighbourhood of $0$. It follows that the curve $\Gamma'$ defined as
$$
\Gamma' = \left\{ \left. \frac{\partial F_t}{\partial t} \right \vert_{t=0} (X,Y,Z) = 0 \right \}
$$
passes through $p_i$ with multiplicity $m_i$ for every $i \ne k$ and through $p_k$ with multiplicity 
$m_k-1$. Indeed, if
$$
\Gamma = \{ f_{m_k}(x,y) + \emph{higher} = 0 \}
$$
in local affine coordinates $(x,y)$ centered at $p_k$, then
\begin{equation}\label{curve}
\Gamma' = \left\{ \dot{a}(0) \frac{\partial f_{m_k}}{\partial x}(x,y) +
\dot{b}(0) \frac{\partial f_{m_k}}{\partial y}(x,y) + \emph{higher} = 0 \right\}.
\end{equation}

Now, if $\Gamma$ has at most two tangent directions at $p_k$, then we may assume that coordinates 
have been chosen so that 
$$
f_{m_k} = x^\alpha y^\beta
$$
with $\alpha, \beta \ge 0$ and $\alpha + \beta = m_k$. It follows that 
$$
\Gamma' = \{ \dot{a}(0) \alpha x^{\alpha-1}y^\beta +
\dot{b}(0) \beta x^\alpha y^{\beta-1} + \emph{higher} = 0 \}.
$$

In particular, by choosing $p_k(t)=[a(t),b(t),1]$ such that one of $\dot{a}(0)$ and $\dot{b}(0)$ is $0$ 
and the other is not $0$ we obtain a curve $\Gamma' \ne \Gamma$ (since their multiplicity at $p_k$ is 
different) of degree $d$ (see \cite{Xu1}, proof of Lemma 1) such that $\Gamma'$ and $\Gamma$ have exactly 
$m_k-1$ tangents in common at $p_k$ (counted with multiplicity).

Hence Bezout's theorem yields
$$
d^2 = \Gamma.\Gamma' \ge \sum_{i \ne k} m_i^2 + m_{k} (m_{k}-1) + m_k -1
= \sum_{i=1}^n m_i^2 -1
$$
and 
$$
C^2 = d^2 - \sum_{i=1}^n m_i^2 \ge - 1.
$$
\qed

\noindent
\emph{Proof of Theorem \ref{two}.} By Theorem \ref{one}, we may assume 
that $m_k=3$ and $\Gamma$ has an ordinary singularity at $p_k$. 
Let $S_k$ be the blow-up of $\mathbb{P}^2$ at the $n-1$ points 
$\{ p_1, \ldots, p_n \} \setminus \{ p_k \}$
and let $\sigma_k: S \to S_k$ be the blow-up of $p_k$. 
If $C_k \subset S_k$ is the strict transform of $\Gamma \subset \mathbb{P}^2$, 
then it is enough to show that $C_k^2 \ge m_k^2 -1$. 

In order to do so, we follow the proof of Lemma 3 in \cite{B} (see also 
\cite{K}, Theorem A). Indeed, the argument of \cite{EL}, Lemma 1.1, and 
of \cite{Xu1}, Lemma 1, as recalled above at the beginning of the proof 
of Theorem \ref{one}, yields non-zero sections $s \in H^0(C,L)$, where 
$$
L := \left( \sigma_k^*(\mathcal{O}_{C_k}(C_k)) \otimes \mathcal{O}_S((1-m_k) E_k) \right) \vert_C
= \mathcal{O}_S(C + E_k) \vert_C.  
$$
In our notation, the sections $s$ are induced by the strict transforms on $S$ of the curves (\ref{curve}). 
Since $\Gamma$ has at least two tangent directions at $p_k$, then as in \cite{Xu2}, proof of Lemma 1, 
Case (1) on p. 202, we have that $\frac{\partial f_{m_k}}{\partial x}(x,y)$ and 
$\frac{\partial f_{m_k}}{\partial y}(x,y)$ are linearly independent modulo higher
degree terms. It follows that the strict transforms on $S$ of 
\begin{eqnarray*}
\Gamma_1' &=& \left\{ \frac{\partial f_{m_k}}{\partial x}(x,y) + \emph{higher} = 0 \right\}, \\
\Gamma_2' &=& \left\{ \frac{\partial f_{m_k}}{\partial y}(x,y) + \emph{higher} = 0 \right\}
\end{eqnarray*}
together with $C + E_k$ generate a net of curves in $\mathbb{P}H^0(S,\mathcal{O}_S(C+E_k))$ 
and induce two linearly independent sections $s_1, s_2 \in H^0(C,L)$. Hence we get a map 
$\phi: C \to \mathbb{P}^1$ of degree $\deg(\phi) \le \deg(L) = C_k^2 - m_k(m_k-1)$. 

Now, if $C$ is rational then it is well-known that $C^2 \ge -1$ (see for instance 
\cite{dF1}, Proposition 2.4). Otherwise, we have
$$
2 \le \deg(\phi) \le C_k^2 - m_k(m_k-1),
$$  
hence
$$
C_k^2 \ge m_k^2 - m_k + 2 = m_k^2 - 1
$$
since $m_k = 3$. 

\qed

\hspace{0.5cm}

\noindent
Claudio Fontanari \newline
Dipartimento di Matematica \newline 
Universit\`a di Trento \newline 
Via Sommarive 14 \newline 
38123 Trento, Italy. \newline
E-mail address: fontanar@science.unitn.it

\end{document}